\documentclass[11pt,reqno]{amsart}

\usepackage{amsfonts,latexsym}
\usepackage{latexsym}
\usepackage{textcomp}
\usepackage{amsmath}
\usepackage{amsfonts}
\usepackage{amssymb}
\usepackage{mathrsfs}

\renewcommand{\a }{\alpha }
\renewcommand{\b }{\beta }
\renewcommand{\d}{\delta }

\newcommand{\D }{\Delta }

\newcommand{\e }{\varepsilon }

\renewcommand{\l }{\lambda }

\newcommand{\m }{\mu }
\newcommand{\n }{\nabla }

\newcommand{\s }{\sigma }
\newcommand{\Sig }{\Sigma}
\renewcommand{\t }{\tau }

\renewcommand{\O }{\Omega }

\newcommand{\ov}{\overline}

\newcommand{\be}{\begin{equation}}
\newcommand{\ee}{\end{equation}}

\newcommand{\R}{\mathbb{R}}
\newcommand{\N}{\mathbb{N}}

\newcommand{\de}{\partial}

\newcommand{\ti}{\tilde}

\newcommand{\M}{\mathcal{M}}

\newcommand{\ra}{{\rangle}}
\newcommand{\la}{{\langle}}

\newcommand{\calO }{\mathcal{O}}

\newcommand{\calL }{\mathcal{L}}

\newcommand{\calU}{{\mathcal U}}

\newtheorem{Theorem}{Theorem}[section]
\newtheorem{Lemma}[Theorem]{Lemma}

\newtheorem{Corollary}[Theorem]{Corollary}

\def\proof{\noindent{{\bf Proof. }}}
\def\square{\vbox{
    \hrule height .4pt
    \hbox{\vrule width .4pt height 7pt \kern 7pt
       \vrule width .4pt}
    \hrule height .4pt }}

\def\QED{\hfill {$\square$}\goodbreak \medskip}

\def\R{{\mathbb R}}

\def\div{{\rm div}}

\font\sc=cmcsc9  \textwidth=14truecm
\title{Weighted Hardy inequality with higher dimensional  singularity
on the boundary }

\author{Mouhamed Moustapha Fall}
\address{\noindent M.M. Fall - Goethe-Universit\"{a}t Frankfurt, Institut
f\"{u}r Mathematik. Robert-Mayer-Str. 10 D-60054 Frankfurt, Germany.} \email{fall@math.uni-frankfurt.de}

\author{Fethi Mahmoudi}
\address{\noindent F. Mahmoudi - Departamento de  Ingenieria Matem\'atica
and CMM, Universidad de Chile, Casilla 170 Correo 3, Santiago,
Chile.}\email{fmahmoudi@dim.uchile.cl}

\begin{document}

\date{}
\maketitle

\bigskip

\noindent {\footnotesize{\bf Abstract.}  Let $\O$ be a smooth
bounded domain in $\R^N$ with   $N\ge 3$ and let $\Sigma_k$ be a
closed smooth submanifold of $\de \O$ of dimension $1\le k\le N-2$.
In this paper we study the weighted Hardy inequality with weight
function singular on $\Sig_k$. In particular we provide 
necessary and sufficient conditions for existence of minimizers.}
\bigskip

\section{Introduction}\label{s:i}

Let $\O$ be a smooth bounded domain of $\R^N$, $N\geq2$ and let
$\Sig_k$ be a smooth closed submanifold of $\de\O$ with dimension
$0\leq k\leq N-1$. Here $\Sig_0$ is a single point and
$\Sig_{N-1}=\de\O$. For $\l\in\R$, consider the problem  of finding
minimizers for the quotient:
\begin{equation}
\label{eq:mpqek} \m_{\l}(\Omega,\Sigma_k):= \inf_{u\in
H^{1}_{0}(\O)} ~\frac{\displaystyle\int_{\O}|\nabla
u|^2p~dx-\l\int_{\O}\d^{-2}|u|^2\eta~dx}
{\displaystyle\int_{\O}\d^{-2}|u|^2q~dx}~,
\end{equation}
where $\d(x):= \textrm{dist}(x,\Sig_k)$ is the distance function to
$\Sig_k$ and where the weights $p,q$ and $\eta$ satisfy
\be\label{eq:weight} \textrm{$p,q\in C^2(\ov{\O})$,}\qquad
p,q>0\quad\textrm{ in $\ov{\O}$,}\qquad \eta>0\quad\textrm{ in
$\ov{\O}\setminus\Sig_k$,}\qquad \textrm{ $\eta\in Lip(\ov{\O})$}
\ee
 and
\begin{equation}\label{eq:min-pq}
\max_{\Sig_k}\frac{q}{p}=1,\qquad  \textrm{ $\eta=0$}\qquad \textrm{ on $\Sig_k$ }.
\end{equation}
 We put
\begin{equation}\label{eq:defIk}
I_{k}=\int_{\Sig_k}\frac{d\s}{\sqrt{1-\left(q(\s)/p(\s)\right)}},\quad
1\leq k\leq N-1\quad\textrm{ and }\quad I_0=\infty.
\end{equation}
It was shown  by Brezis and Marcus in \cite{BM} that there exists $\l^*$ such that
 if $\l>\l^*$ then
$\m_{\l}(\Omega,\Sigma_{N-1})   <\frac{1}{4}$ and it is attained  while for
$\l\leq\l^*$, $\m_{\l}(\Omega,\Sigma_{N-1}) =\frac{1}{4}$ and  it is  not
achieved for every $\l<\l^*$. The critical case
${\l=\l^*}$ was studied by Brezis, Marcus and Shafrir in \cite{BMS}, where they
proved that  $ \m_{\l^*}(\Omega,\Sigma_{N-1})$ admits a minimizer  if and only if
$I_{N-1}<\infty$. The case where $k=0$  ($\Sig_0$ is reduced to a point
 on the boundary) was treated  by  the first author in \cite{Fallccm}
 and the same conclusions hold true.\\
Here we obtain the following
\begin{Theorem}\label{th:mulpqe} Let $\O$ be a smooth bounded
domain of $\R^N$, $N\geq3$  and let $\Sig_k\subset\de\O$ be a closed
submanifold of dimension $k\in[1,N-2]$. Assume that the weight
functions $p,q$ and $\eta$ satisfy \eqref{eq:weight} and
\eqref{eq:min-pq}. Then, there exists
$\l^*=\l^*(p,q,\eta,\O,\Sig_k)$ such that
$$
\begin{array}{ll}
\displaystyle \m_{\l}(\Omega,\Sigma_k)=\frac{(N-k)^2}{4},\quad\forall\l\leq \l^*,\\
\displaystyle \m_{\l}(\Omega,\Sigma_k)<\frac{(N-k)^2}{4},\quad\forall\l> \l^*.
\end{array}
$$
The infinimum   $\m_{\l}(\Omega,\Sigma_k)$ is attained  if $\l>\l^*$ and it is not attained when $\l< \l^*$.
\end{Theorem}
Concerning the critical case we get
\begin{Theorem}\label{th:crit}
Let $\l^*$ be  given by Theorem \ref{th:mulpqe} and consider $I_k$ defined in  \eqref{eq:defIk}. Then
$\m_{\l^*}(\Omega,\Sigma_k)$ is achieved if and only if $I_{k}<\infty $.
\end{Theorem}
By choosing $p=q\equiv1$ and $\eta=\d^2$, we obtain the following consequence of the above theorems.
\begin{Corollary}
Let $\O$ be a smooth bounded domain of $\R^N$, $N\geq3$  and
$\Sig_k\subset\de\O$ be a closed submanifold of dimension
$k\in\{1,\cdots,N-2\}$. For $\l\in\R$, put
$$
\nu_\l(\O,\Sig_k)=\inf_{u\in
H^{1}_{0}(\O)} ~\frac{\displaystyle\int_{\O}|\nabla
u|^2~dx-\l\int_{\O}|u|^2~dx}
{\displaystyle\int_{\O}\d^{-2}|u|^2~dx}~,
$$
Then, there exists $\bar{\l}=\bar{\l}(\O,\Sig_k)$ such that
$$
\begin{array}{ll}
\displaystyle \nu_{\l}(\Omega,\Sigma_k)=\frac{(N-k)^2}{4},\quad\forall\l\leq \bar{\l},\\
\displaystyle \nu_{\l}(\Omega,\Sigma_k)<\frac{(N-k)^2}{4},\quad\forall\l> \bar{\l}.
\end{array}
$$
Moreover $\nu_{\l}(\Omega,\Sigma_k) $ is attained if and only if $ \l> \bar{\l}$.
\end{Corollary}
The proof of the above  theorems are mainly based on the
construction of appropriate sharp $H^1$-subsolution and  $H^1$-supersolutions for the
corresponding operator
$$\calL_\l:=-\D
-\frac{(N-k)^2}{4}q\d^{-2}+\l\d^{-2}\eta $$
 (with $p\equiv 1$).
These super-sub-solutions  are perturbations of an approximate
``virtual" ground-state for the Hardy constant $ \frac{(N-k)^2}{4}$
near $\Sig_k$. For that  we will  consider  the \textit{projection
distance} function $\ti{\d}$ defined near $\Sig_k$ as
$$
\tilde \d(x):=\sqrt{|\mbox{dist}^{\de\O}(\ov
x,\Sigma_k)|^2+|x-\ov x|^2},
$$
where $\ov x$ is  the orthogonal projection of $x$ on $\de\O$ and $\rm{dist}^{\de\O}(\cdot,\Sig_k)$
is the geodesic distance to $\Sig_k$ on $\de\O$  endowed with  the induced metric.
While the distances $\d$ and $\tilde{\d}$  are equivalent, $\D\d$ and $\D\tilde{\d}$
differ and $\d$ does not, in general, provide the right approximate solution for $k\leq N-2$.
Letting $d_{\de\O}=\textrm{dist}(\cdot,\de\O)$, we have
$$
\tilde \d(x):=\sqrt{|\mbox{dist}^{\de\O}(\ov
x,\Sigma_k)|^2+d_{\de \O}(x)^2}.
$$
Our approximate virtual ground-state near $\Sig_k$ reads then as
\be\label{eq:virtgs} x\mapsto d_{\de\O}(x)\,\tilde \d^{
\frac{k-N}{2}}(x). \ee In some appropriate Fermi coordinates
${y}=(y^1,y^2,\dots, y^{N-k}, y^{N-k+1},\dots, y^N)=(\tilde{y},\bar{y})\in\R^{N}$ with $\ti{y}=(y^1,y^2,\dots, y^{N-k})\in\R^{N-k}$ (see next section for
precise definition), the function in \eqref{eq:virtgs} then  becomes
$$
{y}\mapsto y^1|\ti{y}|^{\frac{k-N}{2}}
$$
which is the "virtual" ground-state for the Hardy constant $ \frac{(N-k)^2}{4}$
 in the flat case $\Sig_k= \R^k$ and $\O= \R^N$. We refer to Section \ref{s:pn} for more details about  the constructions of the super-sub-solutions.\\
The proof of the existence part in {Theorem} \ref{th:crit} is inspired from \cite{BMS}. It amounts to obtain a uniform control
of  a specific minimizing sequence for $ \m_{\l^*}(\Omega,\Sigma_k) $ near $\Sig_k$ via the  $H^1$-super-solution constructed.\\
We mention that the existence and non-existence of extremals for
\eqref{eq:mpqek} and related problems were studied in
\cite{AS,CaMuPRSE,CaMuUMI,C,Fall,FaMu,FaMu1,NaC,Na,PT} and some
references therein. We would like to mention that some of the  results in this paper might
of interest in the study of semilinear equations with a Hardy  potential singular
at a submanifold of the boundary. We refer to \cite{Fall-ne-sl, BMR1, BMR2}, 
where existence and nonexistence for semilinear problems were studied via the method
of super/sub-solutions.


%
%
\section{Preliminaries  and Notations}\label{s:pn}
In this section we collect some notations and conventions we are
going to use throughout the paper.

Let $\calU$ be an open subset of $\R^N$, $N\geq 3$, with boundary
$\M:=\de\calU$ a smooth closed hypersurface of ${\R^N}$. Assume that
$\M$ contains a smooth closed submanifold $\Sigma_k$ of dimension
$1\le k\le N-2$. In the following, for $x\in\R^N$, we let $d(x)$ be
the distance function of $\M$ and  $\delta (x)$ the distance
function of $\Sigma_k$.
We denote by $N_\M$ the unit normal vector field of $\M$ pointed into $\calU$.\\
Given $P\in\Sig_k$, the tangent
space $T_P \M$ of $\M$ at $P$ splits as
$$
T_P \M=T_P \Sigma_k\oplus N_P \Sigma_k,
$$
where $T_P\Sigma_k$ is the tangent space of $\Sigma_k$ and $N_P\Sigma_k$ stands for the normal space of $T_P\Sigma_k$ at $P$.
 We assume that the basis of these subspaces are spanned respectively by $\big(E_a\big)_{a=N-k+1,\cdots,N}$ and  $\big(E_i\big)_{i=2,\cdots,N-k} $.
We  will  assume  that $N_\M(P)=E_1$.

A neighborhood of $P$ in $\Sig_k$ can be parameterized via the map
$$
\bar y\mapsto f^P(\bar y)=\textrm{Exp}^{\Sigma_k}_P( \sum_{a=N-k+1}^{N}y^a E_a),
$$
where, $\bar{y}=(y^{N-k+1},\cdots,y^N)$ and where $\textrm{Exp}_P^{\Sigma_k}$
is   the exponential map  at
$P$ in $\Sigma_k$  endowed with
the metric induced by $\M$. Next we extend  $(E_i)_{i=2,\cdots,N-k}$ to an orthonormal frame $(X_i)_{i=2,\cdots,N-k}$ in a neighborhood of $P$.
We can therefore  define  the parameterization of a neighborhood  of $P$  in $\M$ via the mapping
$$
(\breve{y},\bar y)\mapsto h^P_{\M}(\breve{y},\bar y):=\textrm{Exp}^{\M}_{f^P(\bar
y)}\left(\sum_{i=2}^{N-k} y^iX_i\right),
$$
with $ \breve{y}=(y^{2},\cdots,y^{N-k})$
 and   $\textrm{Exp}_Q^\M$ is  the exponential map at $Q$ in $\M$ endowed with
the metric induced by  $\R^N$.
We now have a parameterization of a neighborhood  of $P$  in $\R^N$ defined via the above {Fermi coordinates} by the map
$$
 y=(y^1,\breve{y},\bar y)\mapsto  F^P_{\M}(y^1,\breve{y},\bar y)=h^P_{\M}(\breve{y},\bar y)+y^1 N_\M(h^P_{\M}(\breve{y},\bar y)).
$$
Next we denote by $g$ the metric induced by $F^P_{\M} $ whose components are defined by 
$$g_{\a\b}(y)=\la\de_\a F^P_{\M}(y),\de_\b F^P_{\M}(y)\ra.$$
Then we have the  following expansions (see for instance \cite{FaMah})
\be\label{eq:metexp}
\begin{array}{lll}
g_{11}(y)=1\\
g_{1\b}(y)=0,\quad\quad\quad\quad\quad\quad\textrm{ for } \b=2,\cdots,N\\
g_{\a\b}(y)=\d_{\a\b}+\calO(|\tilde{y}|),\quad\textrm{ for } \a,\b=2,\cdots,N,
\end{array}
\ee
where $\tilde{y}=(y^1,\breve{y})$ and $\calO(r^m)$ is a smooth function in the variable $y$ which is uniformly bounded by
a constant (depending only $\M$ and $\Sig_k$) times  $r^m$.

In concordance to the above coordinates, we will consider the ``half"-geodesic neighborhood contained in $\calU$ around
$\Sigma_k$ of radius $\rho$
\be\label{eq:geodtub}
\calU_{\rho}(\Sigma_k) := \{ x \in \calU: \quad \ti{\d}(x)<\rho \},
\ee
with $\tilde \d $ is the projection distance function given by
$$
\tilde \d(x):=\sqrt{|\mbox{dist}^{\M}(\ov
x,\Sigma_k)|^2+|x-\ov x|^2},
$$
where $\ov x$ is  the orthogonal projection of $x$ on $\M$ and $\rm{dist}^{\M}(\cdot,\Sig_k)$
is the geodesic distance to $\Sig_k$ on $\M$ with the induced metric.
Observe that
\be\label{eq:tidFptiy}
\tilde \d(F^P_\M(y))=|\tilde y|,
\ee
where $\tilde y=(y^1,\breve{y})$.
 We also
define $\sigma(\ov x)$ to be the orthogonal projection of $\ov x$ on $\Sigma_k$ within $\M$.
Letting
$$
\hat \delta(\ov x):=\mbox{dist}^{\M}(\ov x,\Sigma_k),
$$
 one has
$$
\ov x=\textrm{Exp}_{\sigma(\ov x)}^\M(\hat\d\,\n\hat\d)\quad \hbox{or
equivalently }\quad \sigma(\ov x)=\textrm{Exp}_{\ov x}^\M(-\hat\d\,\n\hat\d).
$$
Next we observe that
\be\label{eq:td-hd}
\ti{\d}(x)=\sqrt{\hat{\d}^2(\bar{x})+d^2(x)}.
\ee
 In addition it can be easily checked via the implicit function theorem that there exists a positive constant
 $\b_0=\b_0(\Sig_k,\O)$ such that $\ti{\d}\in C^\infty(\calU_{\b_0}(\Sig_k))$.

 It is clear  that for
$\rho$ sufficiently small, there exists a finite number of Lipschitz
open sets $(T_i)_{1\le i\le N_0}$ such that
$$
T_i\cap T_j=\emptyset \quad \hbox{for }\,i\ne j\quad \hbox{and}\quad
\calU_\rho(\Sig_k)=\bigcup_{i=1}^{N_0}\ov{ T_i}.
$$
We may assume that each $T_i$ is chosen, using the above coordinates, so that
$$
T_i=F^{p_i}_{\M}(B^{N-k}_+(0,\rho)\times D_i)\quad\hbox{with }\; p_i\in \Sigma_k,
$$
where the $D_i$'s are Lipschitz disjoint open sets of $\R^k$ such that
$$
\bigcup_{i=1}^{N_0} \ov{f^{p_i} (D_i)}=\Sig_k.
$$
In the above setting we have
\begin{Lemma} \label{lemddelta} As $\ti{\d}\to0$, the following expansions hold
\begin{enumerate}
\item $\d^2=\ti{\d}^2(1+O(\ti{\d}))$,
\item $\n \ti{\d}\cdot\n d=\displaystyle\frac{d}{\ti{\d}}$,
\item $|\n\tilde{\d}|=1+O(\ti{\d}),$
\item $\Delta \ti{\delta }=\frac{N-k-1}{\ti{\delta}}+O(1)$,
\end{enumerate}
where $O(r^m)$ is a function  for which there exists a constant $C=C(\M,\Sig_k)$ such that
$$
|O(r^m)|\leq C r^m.
$$
\end{Lemma}

\proof
\begin{enumerate}
\item Let $P\in \Sig_k$. With an abuse of notation, we write $x(y)= F^P_\M(y)$ and we  set
$$
\vartheta( y):=\frac12\delta^2 (x({y})).
$$
 The function $\vartheta$ is
smooth in a small neighborhood of the origin in $\R^{N}$ and Taylor
expansion yields
\begin{eqnarray}
  \vartheta( y)&=&\vartheta(0,\bar{y})\tilde y+\nabla\vartheta(0,\bar{y})[\tilde y]+\frac12\nabla^2\vartheta(0,\bar{y})[\tilde y,\tilde y]+\calO(\|\tilde y\|^3)\nonumber\\
&=&\label{eq:vartzyb}\frac12\nabla^2\vartheta(0,\bar{y})[\tilde y,\tilde y]+\calO(\|\tilde y\|^3) .
\end{eqnarray}
Here we have used the fact that $x(0,\bar{y} )\in \Sig_k$ so that $ \d(x(0,\bar{y}))=0$.
We write
$$
\nabla^2\vartheta(0,\bar{y})[\tilde y,\tilde
y]=\sum_{i,l=1}^{N-k}\Lambda_{il}y^iy^l,
$$
with
\begin{eqnarray*}
      \Lambda_{il} &:=&\frac{\partial^2 \vartheta}{\partial y^i\partial y^l}/_{ \ti{y}=0}\\
       &=& \frac{\partial}{\partial y^l}\bigg(\frac{\partial }{\partial x^j} (\frac12 \delta^2(x)\frac{\partial x^j}{\partial y^i} ) \bigg)/_{
       \ti{y}=0}\\
       &=&\frac{\partial^2}{\partial x^i\partial x^s}(\frac12
       \delta^2)(x)\frac{\partial {x^j}}{\partial y^i}\frac{\partial x^s}{\partial y^l}/_{
       \ti{y}=0}+\frac{\partial }{\partial x^j}(\delta^2)(x)\frac{\partial^2x^s}{\partial y^i\partial y^l}/_{
       \ti{y}=0}.
    \end{eqnarray*}
Now using the fact that
$$
\frac{\partial x^s}{\partial y^l}/_{ \ti{y}=0}=g_{ls}=\delta_{ls}\quad
\textrm{and}\quad\frac{\partial }{\partial x^j}(\delta^2)(x)/_{
       \ti{y}=0}=0,
$$
we obtain
\begin{eqnarray*}
      \Lambda_{il} y^i y^l&=&y^i y^s\,\frac{\partial^2}{\partial x^i\partial x^s}(\frac12
       \delta^2)(x)/_{ \ti{y}=0} \\
       &=& |\tilde y|^2,
    \end{eqnarray*}
where we have used the fact that the matrix $\left(\frac{\partial^2}{\partial x^i\partial x^s}(\frac12
       \delta^2)(x)/_{ \ti{y}=0} \right)_{1\leq i,s\leq N}$ is the matrix of the orthogonal projection  onto the normal space of $T_{f^P(\bar{y})}\Sig_k$.
Hence using \eqref{eq:vartzyb}, we get
$$
\delta^2 (x({y}))=|\tilde y|^2 +\calO(|\tilde y|^3).
$$
This together with \eqref{eq:tidFptiy} prove the first expansion.

 \item Thanks to \eqref{eq:tidFptiy} and \eqref{eq:metexp}, we infer that
$$
\n \ti{\d}\cdot\n d(x(y))= \frac{\de \ti{\d}( x(y))}{\de y^1}=\frac{y^1}{|\ti{y}|}=\frac{d(x(y))}{\ti{\d}(x(y))}
$$
as desired.
\item We observe that
$$
\frac{\de \ti{\d}}{\de x^\t}\frac{ \de \ti{\d}}{\de x^\t} (x(y)) =g^{\t \a}(y)g^{\t \b}(y)\frac{\de \ti{\d}(x(y))}{\de y^\a}\frac{\de \ti{\d}(x(y))}{\de y^\b},
$$
where $(g^{\a\b})_{\a,\b=1,\dots,N} $ is  the inverse of the  matrix $(g_{\a\b})_{\a,\b=1,\dots,N} $.
Therefore using  \eqref{eq:tidFptiy} and \eqref{eq:metexp}, we get the result.

\item Finally using  the expansion of the  Laplace-Beltrami operator  $\D_g$, see Lemma 3.3 in  \cite{mm}, applied to \eqref{eq:tidFptiy}, we get
 the last estimate.
\end{enumerate}
\QED
In the following of -- only -- this section, $q:\ov{\calU} \to \R$
be such that \be\label{eq:q} q\in C^2(\ov{\calU}),\quad\textrm{ and
}\quad q\leq 1\quad\textrm{ on } \Sig_k. \ee Let $M,a\in\R$, we
consider  the function \be\label{eq:pert-gst}
W_{a,M,q}(x)=X_a(\ti{\delta}(x))\,e^{Md(x)}\,d(x)\,\ti{\delta}(x)^{\alpha(x)},
\ee where

$$
X_a(t)=(-\log(t))^a\quad \, 0<t<1 $$
and 

$$
\alpha(x)=\frac{k-N}{2}+\frac{N-k}{2}\sqrt{1-q(\s(\bar{x}))+\ti{\d}(x)}.
$$

In the above setting,  the following useful result holds.
\begin{Lemma}\label{LapFinalExp}
 As $\d\to 0$,   we have
\begin{eqnarray*}
  \Delta W_{a,M,q}&=& - \frac{(N-k)^2}{4}\,q\,\delta^{-2} \,W_{a,M,q}-{2\,a\,\sqrt{\tilde\alpha}}\,X_{-1}(\d)\,\delta^{-2}\,W_{a,M,q}
 \\
  &+& {a(a-1)} \,X_{-2}(\d)\,\delta^{-2}\,W_{a,M,q}+\frac{h+2M}{d}\,W_{a,M,q}+O(|\log(\delta)|\,\delta^{-\frac32})\,W_{a,M,q},\nonumber
\end{eqnarray*}
where $\tilde{\alpha}(x)=\frac{(N-k)^2}{4}\left(1- q(\sigma (\ov x))+\ti{\delta} (x)\right) $ and $h=\D d$. Here the lower order term satisfies
$$
|O(r)|\leq C |r|,
$$
where $C$ is a positive constant  only depending  on $a,M,\Sig_k,\calU$ and $\|q\|_{C^2(\calU)}$.
\end{Lemma}
\proof
We put $s=\frac{(N-k)^2}{4} $.
Let $w=\ti{\delta}(x)^{\alpha(x)} $ then the following formula can be easily verified
\be\label{eq:1}
\D w=w\bigg( \D \log(w)+|\nabla\log(w)|^2 \bigg).
\ee
Since
$$
\log(w)=\alpha\log(\ti{\delta}),
$$
we get
\be\label{eq:2}
\D \log(w)=\D \alpha\log(\ti{\delta})+2\nabla\alpha\cdot \nabla
(\log(\ti{\delta}))+\alpha\D \log(\ti{\delta}).
\ee
We have
\be\label{eq:3}
\D\alpha=\D\sqrt{\tilde \alpha}=\sqrt{\tilde
\alpha}\,\left(\frac12 \D\log(\tilde \alpha) +\frac14|\nabla
\log(\tilde \alpha)|^2 \right),
\ee
$$
\nabla\log(\tilde\alpha)=\frac{\nabla\tilde\alpha}{\tilde\alpha}=\frac{-s\nabla(q\circ\sigma)+s\nabla\ti{\delta}}{\tilde\a}
$$
and using the formula \eqref{eq:1}, we obtain
\begin{eqnarray*}
  \D\log(\tilde\alpha)&=&\frac{\D\tilde\alpha}{\tilde\alpha} -\frac{|\nabla\tilde\alpha|^2}{\tilde\alpha^2}\\
   &=& \frac{-s\D(q\circ\sigma)+s\D\ti{\delta}}{\tilde\alpha} -\frac{s^2|\nabla(q\circ\sigma)|^2+s^2|\nabla\ti{\delta}|^2}
{\tilde\alpha^2}+2s^2\frac{\nabla(q\circ\sigma)\cdot\nabla\ti{\delta}}{\tilde\alpha^2}.
\end{eqnarray*}
Putting the above in \eqref{eq:3}, we deduce that
\be\label{eq:4}
  \D\alpha =\frac{1}{2\sqrt{\tilde\alpha}} \bigg( -s\D (q\circ\sigma)+s\D\ti{\delta}-
  \frac12\frac{s^2|\nabla(q\circ\sigma)|^2+s^2|\nabla\ti{\delta}|^2-2s^2\nabla(q\circ\sigma)\cdot\nabla\ti{\delta}}{\tilde\alpha}\bigg).
\ee
Using Lemma \ref{lemddelta} and
the fact that $q$ is in $C^2(\ov{\calU})$,
 together with \eqref{eq:4}  we get
\be\label{eq:5}
\D\alpha= O({\ti{\delta}^{-\frac32}}).
\ee
 On the other hand
$$
\nabla
\alpha=\nabla\sqrt{\tilde\alpha}=\frac12\frac{\nabla\tilde\alpha}{\sqrt{\tilde\alpha}}=-\frac{s}{2\sqrt{\tilde\alpha}}\nabla(q\circ\sigma)+
\frac{s}{2}\frac{\nabla\ti{\delta}}{\sqrt{\tilde\alpha}}
$$
so that
$$
\nabla \alpha\cdot \nabla
\ti{\delta}=-\frac{s}{2\sqrt{\tilde\alpha}}\nabla(q\circ\sigma)\cdot
\nabla \ti{\delta}+
\frac{s}{2}\frac{|\nabla\ti{\delta}|^2}{\sqrt{\tilde\alpha}}=O(\ti{\d}^{-\frac12})
$$
and from which we deduce that
\be\label{eq:6}
  \nabla\alpha\cdot \nabla\log(\ti{\delta}) = \frac{1}{\ti{\delta}} \nabla\alpha\cdot \nabla\ti{\delta}
  =
O(\ti{\d}^{-\frac32}).
\ee
By Lemma \ref{lemddelta} we have that
$$
\alpha\D\log(\ti{\delta})=\alpha\,\frac{N-k-2}{\ti{\delta}^2}\,(1+O(\ti{\delta})).
$$
Taking back the above estimate together with \eqref{eq:6} and \eqref{eq:5} in \eqref{eq:2}, we get
\be\label{eq:7}
  \D\log(w) =  \alpha\,\frac{N-k-2}{\ti{\delta}^2}\,(1+O(\ti{\delta}))
  +O(|\log(\ti{\d})|\ti{\d}^{-\frac32}).
\ee
We also have
$$
\nabla(\log(w))=\nabla(\alpha \log(\ti{\delta}))=\alpha
\frac{\nabla\ti{\delta}}{\ti{\delta}}+\log(\ti{\delta})\nabla \alpha
$$
and thus
$$
|\nabla(\log(w))|^2=\frac{\alpha^2}{\ti{\delta}^2}+\frac{2\alpha\log(\ti{\delta})}{\ti{\delta}}\,\nabla\ti{\delta}\cdot\nabla
\alpha+|\log(\ti{\delta})|^2|\nabla \alpha|^2=\frac{\alpha^2}{\ti{\delta}^2}+ O(|\log(\ti{\d})|\ti{\d}^{-\frac32}).
$$
Putting this together with \eqref{eq:7} in \eqref{eq:1}, we conclude that
\be\label{eq:8}
 \frac{ \D w }{w}=
 \alpha\,\frac{N-k-2}{\ti{\delta}^2}+\frac{\alpha^2}{\ti{\delta}^2}+O(|\log(\ti{\delta})|\,\ti{\delta}^{-\frac32}).
\ee
Now we define the function
$$
v(x):=d(x)\,w(x),
$$
where we recall that $d$ is the distance function to the boundary of $\calU$.
It is clear that
\be\label{eq:9}
\D v= w\D d+d\D w+2\nabla d\cdot \nabla w.
\ee
Notice that
$$
\nabla w=w\,\nabla
\log(w)=w\,\left(\log(\ti{\delta})\nabla\alpha+\alpha\frac{\nabla
\ti{\delta}}{\ti{\delta}}\right)
$$
and so
\be\label{eq:10}
\nabla d\cdot\nabla w=w\,\left(\log(\ti{\delta})\nabla d
\cdot\nabla\alpha+\frac{\alpha}{\ti{\delta}}\nabla d\cdot\nabla
\ti{\delta}\right).
\ee
Recall the second assertion  of  Lemma \ref{lemddelta} that we rewrite as 
\be\label{eq:11}
\nabla d\cdot\nabla \ti{\delta}=\frac{d}{\ti{\delta}}.
\ee
Therefore
\be\label{eq:12}
\nabla d \cdot\nabla\alpha=\nabla
d\cdot\left(-\frac{s}{2\sqrt{\tilde
\alpha}}\nabla(q\circ\sigma)+\frac{s}{2}\frac{\nabla\ti{\delta}}{\sqrt{\tilde
\alpha}} \right)=\frac{s}{2\sqrt{\tilde
\alpha}}\frac{d}{\ti{\delta}}-\frac{s}{2\sqrt{\tilde
\alpha}}\nabla d\cdot\nabla(q\circ\sigma).
\ee
Notice that if $x$ is in a neighborhood of some point $P\in \Sig_k$ one has
$$
\nabla d\cdot\nabla(q\circ\sigma)(x)=\frac{\de}{\de
y^1}q(\sigma(\ov x))=\frac{\de}{\de y^1}q( f^P(\ov y))=0.
$$
This with \eqref{eq:12} and \eqref{eq:11} in \eqref{eq:10} give
\begin{eqnarray}\label{eq:13}
  \nabla d\cdot\nabla w&=&w\,\left(O(\ti{\delta}^{-\frac32}|\log(\ti{\delta})|)\,d+\frac{\alpha}{\ti{\delta}^2}\,
  d \right)\nonumber \\
  &=& v\,\left(O(\ti{\delta}^{-\frac32}|\log(\ti{\delta})|)+\frac{\alpha}{\ti{\delta}^2}\right).
\end{eqnarray}
From \eqref{eq:8}, \eqref{eq:9} and \eqref{eq:13} (recalling the expression of $\a$ above), we  get immediately
\begin{eqnarray}\label{eq:14}
  \D v&=&\left(
\alpha\,\frac{N-k}{\ti{\delta}^2}+\frac{\alpha^2}{\ti{\delta}^2}\right)\,v+O(|\log(\ti{\delta})|\,\ti{\delta}^{-\frac32})\,v+
\frac{h}{d}\,v \nonumber\\
&=&\left(- \frac{(N-k)^2}{4}
\frac{q(x)}{\ti{\delta}^2}+O(|\log(\ti{\delta})|\,\ti{\delta}^{-\frac32})\right)\,v+
\frac{h}{d}\,v,
\end{eqnarray}
where $h=\D d$. Here we have used the fact that $|q(x)-q(\s(\bar{x}))|\leq C \ti{\d }(x)$  for $x$ in a neighborhood of $\Sig_k$.\\
Recall that
$$
W_{a,M,q}(x)=X_a(\ti{\delta}(x))\,e^{Md(x)}\,v(x), \quad \hbox{ with }\quad
X_a(\ti{\delta}(x)):=(-\log(\ti{\delta}(x)))^a,
$$
where $M$ and $a$ are two real numbers. We have
\begin{eqnarray*}
  \D W_{a,M,q} = X_a(\ti{\delta})\,\D (e^{Md}\,v)+2\nabla X_a(\ti{\delta})\cdot\nabla (e^{Md}\,v)+e^{Md}\,v\,\D X_a(\ti{\delta})
\end{eqnarray*}
and thus
\be\label{eq:15}
\begin{array}{lll}
  \D W_{a,M,q}
  &= &X_a(\ti{\delta})e^{Md}\,\D
  v+X_a(\ti{\delta}) \D (e^{Md})\, v+2X_a(\ti{\delta})\n v\cdot \nabla(e^{Md})\\
  &\,\,&+\,2\nabla X_a(\ti{\delta})\cdot\left( v\,\nabla (e^{Md})+e^{Md}\nabla v\right)+e^{Md}\,v\,\D
  X_a(\ti{\delta}).
\end{array}
\ee
We shall estimate term by term the above expression.\\
 First we have form \eqref{eq:14}
\be\label{eq:141}
X_a(\ti{\delta})e^{Md}\,\D v= - \frac{(N-k)^2}{4}
\frac{q}{\ti{\delta}^2}\, W_{a,M,q} +
\frac{h}{d}\, W_{a,M,q} +O(|\log(\ti{\delta})|\,\ti{\delta}^{-\frac32})\, W_{a,M,q}.
\ee
It is plain  that
\be\label{eq:17}
X_a(\ti{\delta})\,\D (e^{Md})\,v=O(1)
\,W_{a,M,q}.
\ee
It is clear  that
\be\label{eq:nv}
\nabla v= w\,\nabla d+d\,\nabla w=w\,\nabla
d+d\,\left(\log(\ti{\delta})\,\nabla\alpha+\alpha \frac{\nabla
\ti{\delta}}{\ti{\delta}}\right)\, w.
\ee
From which and \eqref{eq:11} we get
\begin{eqnarray}\label{eq:16}
  X_a(\ti{\delta})\,\nabla v\cdot \nabla(e^{Md}) &=& M\,X_a(\ti{\delta})\,e^{Md}\,w\left\{ |\nabla d|^2+d\, \left(\log(\ti{\delta})\,\nabla d\cdot \nabla\alpha+
 \frac{\alpha }{\ti{\delta}}
  \nabla\ti{\delta}\cdot\nabla d\right)\right\}\nonumber \\
  &=&M\,X_a(\ti{\delta})\,e^{Md}\,w\left\{
  1+O(|\log(\ti{\delta})|\,\ti{\delta}^{-\frac12})\,d+O(\ti{\delta}^{-1})\,d\right\}\nonumber\\
  &=& W_{a,M,q} \,\left\{ \frac{M}{d}+O(|\log(\ti{\delta})|\,\ti{\delta}^{-1})\right\}.
\end{eqnarray}
Observe that
$$
\nabla(X_a(\ti{\delta}))=-a\,\frac{\nabla \ti{\delta}}{\ti{\delta}}
X_{a-1}(\ti{\delta}).
$$
This with \eqref{eq:nv} and \eqref{eq:11} imply that
\be\label{eq:18}
\nabla X_a(\ti{\delta})\cdot\left( v\,\nabla (e^{Md})+e^{Md}\nabla
v\right)=
-\frac{a(\alpha+1)}{\ti{\delta}^2}\,X_{-1}\,W_{a,M,q}+O(|\log(\ti{\delta})|\ti{\delta}^{-\frac32})\,W_{a,M,q}.
\ee
By Lemma \ref{lemddelta}, we have
$$
\D(X_a(\ti{\delta}))=\frac{a}{\ti{\delta}^2}X_{a-1}(\ti{\delta})\{2+k-N+O(\ti{\delta})\}+\frac{a(a-1)}{\ti{\delta}^2}X_{a-2}(\ti{\delta}).
$$
Therefore we obtain
\be\label{eq:19}
e^{Md}v \D(X_a(\ti{\delta}))=\frac{a}{\ti{\delta}^2} \{2+k-N+O(\ti{\delta})\}\,X_{-1}\,W_{a,M,q}+ \frac{a(a-1)}{\ti{\delta}^2}X_{-2} \,W_{a,M,q}.
\ee
Collecting  \eqref{eq:141},  \eqref{eq:17}, \eqref{eq:16}, \eqref{eq:18} and \eqref{eq:19} in the expression \eqref{eq:15},
 we get
  as $\ti{\d}\to 0$
\begin{eqnarray*}
  \Delta W_{a,M,q}&=& - \frac{(N-k)^2}{4}\,q\,\ti{\delta}^{-2} \,W_{a,M,q}-2\,a\,\sqrt{\ti{\alpha}}\,X_{-1}(\ti{\d})\,\ti{\delta}^{-2}\,W_{a,M,q}
 \\
  &+& {a(a-1)} \,X_{-2}(\ti{\d})\,\ti{\delta}^{-2}\,W_{a,M,q}+\frac{h+2M}{d}\,W_{a,M,q}
  +O(|\log(\ti{\delta})|\,\ti{\delta}^{-\frac32})\,W_{a,M,q}.\nonumber
\end{eqnarray*}
 The conclusion of the  lemma  follows at once from the first assertion  of Lemma \ref{lemddelta}.
\QED

\subsection{Construction of a  subsolution}
For $\l\in\R$ and  $\eta\in Lip(\ov{\calU})$ with $\eta=0$ on $\Sig_k$, we define the operator
\be\label{eq:calL_l}
\mathcal{L}_\l:=
-\D -\frac{(N-k)^2}{4}\,q\,\delta^{-2}+\l\, \eta\,\delta^{-2},
\ee
where $q$ is as in \eqref{eq:q}.
We have the following lemma
\begin{Lemma} \label{le:lowerbound}
There exist two positive constants $M_0,\beta_0$  such that for all
$\beta\in\,(0,\beta_0)$ the function
$V_\e:=W_{-1,M_0,q}+W_{0,M_0,q-\e}$ (see \eqref{eq:pert-gst})
satisfies
\begin{equation}\label{eq:subsolution}
\mathcal{L}_\l V_\e\le 0 \quad \textrm{ in } \calU_\b,\quad\hbox{ for all }\; \e\in[0,1).
\end{equation}
Moreover $V_\e\in H^1(\calU_\beta)$ for any $\e\in(0,1)$ and in addition
\begin{equation}\label{eq:Iq}
\int_{\calU_\b}\frac{V_{0}^2}{\d^2}\,dx\geq C \int_{\Sigma_k} \frac{1}{\sqrt{1-q(\sigma)}}\,d\sigma.
\end{equation}
\end{Lemma}
\proof Let $\beta_1$ be a positive small real number so that $d$ is
smooth in $\calU_{\b_1}$.  We choose
$$
 M_0= \max_{x\in \ov\calU_{\b_1}}|h(x)|+1.
$$
Using this and Lemma \ref{LapFinalExp}, for some $\b\in(0,\b_1)$, we have
\be\label{eq:LaM0}
  \mathcal{L}_\l W_{-1,M_0,q} \le \left(-2\delta^{-2} \,X_{-2}+C|\log(\delta)|\,\delta^{-\frac32}+|\l|\eta \d^{-2}\right)\,W_{-1,M_0,q}\quad
\textrm{ in } \calU_\b. \ee
Using the
fact that   the function $\eta$ vanishes on
$\Sigma_k$ (this implies in particular that $|\eta|\le C \delta$ in
$\calU_\b$), we have
$$
\mathcal{L}_\l(W_{-1,M_0,q})\le -\delta^{-2} \,X_{-2}\,W_{-1,M_0,q}= -\delta^{-2} \,X_{-3}\,W_{0,M_0,q}\quad \textrm{ in }\calU_\b,
$$
for $\b$ sufficiently small. Again by Lemma \ref{LapFinalExp}, and
similar arguments as above, we have \be\label{eq:LaMqep}
  \mathcal{L}_\l W_{0,M_0,q-\e} \le C|\log(\delta)|\,\delta^{-\frac32}\,W_{0,M_0,q-\e}\leq C|\log(\delta)|\,\delta^{-\frac32}\,W_{0,M_0,q}\quad\textrm{ in }\calU_{\b},
\ee
for any $\e\in [0,1)$. Therefore we get
$$
 \mathcal{L}_\l \left(W_{-1,M_0,q}+W_{0,M_0,q-\e} \right)\leq 0\quad \textrm{ in }\calU_{\b},
$$
if $\b$ is small.  This proves \eqref{eq:subsolution}.\\
The proof of the fact that
$W_{a,M_0,q}\in H^1(\calU_\beta)$, for any $a<-\frac{1}{2}$ and $ W_{0,M_0,q-\e}\in H^1(\calU_\beta) $, for $\e>0$ can be easily checked using polar coordinates
(by assuming without any loss of generality that $M_0=0$ and $q\equiv 1$), we therefore skip it.  \\
We now prove the last statement of the theorem.
Using Lemma \ref{lemddelta}, we have

\begin{eqnarray*}
\int_{\calU_\b}\frac{V_{0}^2}{\d^2}\,dx
&\ge& \int_{\calU_\b}\frac{W_{0,M_0,q}^2}{\d^2}\,dx\\
&\ge &C\,\int_{\calU_\b(\Sig_k)}d^2(x)\tilde{\delta}(x)^{2\a(x)-2}\,dx\\
&\ge& C\sum_{i=1}^{N_0}\,\int_{T_i}d^2(x)\tilde{\delta}(x)^{2\a(x)-2}\,dx\\
&=&C\sum_{i=1}^{N_0}\,\int_{B^{N-k}_+(0,\b)\times
D_i}(y^1)^2\,|\tilde y|^{2\a(F^{p_i}_\M(y))-2}\,
|{\rm Jac}(F^{p_i}_\M)|(y)\,dy\\
&\ge& C\,\sum_{i=1}^{N_0}\,\int_{B^{N-k}_+(0,\b)\times
D_i}(y^1)^2\,|\tilde y|^{k-N-2+(N-k)\sqrt{1-q(f^{p_i}(\bar
y))}}\, \,|\tilde y|^{-\sqrt{|\ti{y}|}}\,dy.
 \end{eqnarray*}
Here we used the fact that $|{\rm Jac}(F^{p_i}_\M)|(y)\ge C$. Observe that
$$
|\tilde y|^{-\sqrt{|\ti{y}|}}\ge C >0
\quad \hbox{as }\, |\tilde y| \to 0.
$$
Using polar coordinates, the above integral becomes
\begin{eqnarray*}
\int_{\calU_\b}\frac{V_{0}^2}{\d^2}\,dx &\ge&
 C\,\sum_{i=1}^{N_0}\int_{D_i}\int_{S^{N-k-1}_+}\left(\frac{y^1}{|\tilde
y|}\right)^2\,d
\theta\int_0^{\b}r^{-1+(N-k)\sqrt{1-q(f^{p_i}(\bar
y))}}\,d\bar y
\\
&\ge & C\,\sum_{i=1}^{N_0}\int_{D_i}\int_0^{r_{i_1}}r^{-1+(N-k)\sqrt{1-q(f^{p_i}(\bar
y))}}\,|\textrm{Jac}(f^{p_i})|(\bar y)\,d\bar y.
 \end{eqnarray*}
We therefore obtain
\begin{eqnarray*}
\int_{\calU_\b}\frac{V_{0}^2}{\d^2}\,dx
&\geq &  C\,\int_{\Sig_k}\int_0^{\b}r^{-1+(N-k)\sqrt{1-q(\s)}}\,dr\,d\s\\
&\geq &  C\,\int_{\Sig_k}\frac{1}{\sqrt{1-q(\s)}}\,d\s.
 \end{eqnarray*}
 This concludes the proof of the lemma.
\QED

\subsection{Construction of a  supersolution}
In this subsection we provide a supersolution for the operator $\calL_\l$ defined in \eqref{eq:calL_l}. We prove
\begin{Lemma} \label{le:upperbound}
 There exist  constants $\beta_0>0$,
 $M_{1}<0,$ $M_0>0$ (the constant $M_0$ is as in Lemma \ref{le:lowerbound}) such that for
all $\beta\in\,(0,\beta_0)$ the function $U:=W_{0,M_1,q}-W_{-1,M_0,q}>0$ in $\calU_\b$ and
satisfies
\begin{equation}\label{eq:supsolution}
\mathcal{L}_\l U_a \geq 0  \quad \textrm{ in } \calU_\b.
\end{equation}
Moreover $U\in H^1(\calU_\beta)$
provided
\begin{equation}\label{eq:Iql}
\int_{\Sigma_k} \frac{1}{\sqrt{1-q(\sigma)}}\,d\sigma <+\infty.
\end{equation}
\end{Lemma}
\proof
 We consider $\b_1$ as in the beginning of the proof of Lemma \ref{le:lowerbound} and we define
\begin{equation}\label{eq:M1}
 M_1=-\frac12\,\max_{x\in\ov\calU_{\beta_1}}|h(x)|-1.
\end{equation}
Since $$ U(x)=(e^{M_1 d(x)}-e^{M_0d(x)}X_{-1}(\ti{\d}(x)))d(x)\ti{\d}(x)^{\a(x)},$$ it follows that $U>0$  in $\calU_\b$ for   $\b>0$ sufficiently small.
By \eqref{eq:M1} and Lemma \ref{LapFinalExp}, we get
\begin{eqnarray*}
  \mathcal{L}_\l W_{0,M_1,q} \ge \left(-C|\log(\delta)|\,\delta^{-\frac32}-|\l|\eta \d^{-2}\right)\,W_{0,M_1,q}.
\end{eqnarray*}
Using  \eqref{eq:LaM0} we have
$$
  \mathcal{L}_\l (- W_{-1,M_0,q})\geq
\left(2\d^{-2}X_{-2}-C|\log(\delta)|\,\delta^{-\frac32}-|\l|\eta \d^{-2}\right)\, W_{-1,M_0,q}.
$$
Taking the sum of the two above inequalities, we obtain
$$
\mathcal{L}_\l U\geq0 \quad\textrm{ in }\calU_\b,
$$
which holds true because $|\eta|\leq C\d$ in $ \calU_{\b}$.
Hence we get readily \eqref{eq:supsolution}.\\
Our next task is to prove  that $U\in H^1(\calU_\b)$ provided \eqref{eq:Iql} holds, to do so it is enough to show that $W_{0,M_1,q} \in H^1(\calU_\b)$ provided \eqref{eq:Iql} holds.\\
We  argue as in the proof of  Lemma \ref{le:lowerbound}. We have
\begin{eqnarray*}
 \int_{\calU_\b}|\nabla W_{0,M_1,q}|^2 &\le & C\int_{\calU_\b}d^2(x)\ti{\delta}(x)^{2\a(x)-2}\,dx\\
&\leq& C\sum_{i=1}^{N_0}\int_{B^{N-k}_+(0,\b)\times
D_i}d^2(F^{p_i}_\M(y))\tilde{\delta}(F^{p_i}_\M(y))^{2\a(F^{p_i}_\M(y))-2}
|{\rm
Jac}(F^{p_i}_\M)|(y)dy\\
&\leq&C\sum_{i=1}^{N_0}\,\int_{B^{N-k}_+(0,\b)\times
D_i}(y^1)^2\,|\tilde y|^{2\a(F^{p_i}_\M(y))-2}\,
|{\rm Jac}(F^{p_i}_\M)|(y)\,dy\\
&\le& C\,\sum_{i=1}^{N_0}\,\int_{B^{N-k}_+(0,\b)\times
D_i}(y^1)^2\,|\tilde y|^{k-N-2+(N-k)\sqrt{1-q(f^{p_i}(\bar
y))}}\, \,|\tilde y|^{-\sqrt{|\ti{y}}|}\,dy.
 \end{eqnarray*}
Here we used the fact that $|{\rm Jac}(F^{p_i}_\M)|(y)\le C$.  Note that
$$
|\tilde y|^{-\sqrt{|\ti{y}}|}\le C
\quad \hbox{as }\, |\tilde y|\to 0.
$$
Using polar coordinates, it follows that
\begin{eqnarray*}
\int_{\calU_\b}|\nabla W_{0,M_1,q}|^2
&\le& C\,\sum_{i=1}^{N_0}\int_{D_i}\int_{S^{N-k-1}_+}\left(\frac{y^1}{|\tilde
y|}\right)^2\,d
\theta\int_0^{\b}r^{-1+(N-k)\sqrt{1-q(f^{p_i}(\bar
y))}}\,dr\,d\bar y\\
&\le&
C\,      \sum_{i=1}^{N_0}\,\int_{D_i}\frac{1}{\sqrt{1-q(f^{p_i}(\bar
y))}}\,d\bar y.
 \end{eqnarray*}
Racalling that $|{\rm Jac}(f^{p_i})|(\bar y)=1+O(|\bar y|)$, we deduce that
\begin{eqnarray*}
\sum_{i=1}^{N_0}\,\int_{D_i}\frac{1}{\sqrt{1-q(f^{p_i}(\bar
y))}}\,d\bar y&\le&
C\sum_{i=1}^{N_0}\,\int_{D_i}\frac{1}{\sqrt{1-q(f^{p_i}(\bar
y))}}\,|{\rm Jac}(f)|(\bar y)\,d\bar
y\\
&=&C\int_{\Sigma_k}\frac{1}{\sqrt{1-q(\sigma})}\,d\sigma.
 \end{eqnarray*}
Therefore
\begin{eqnarray*}
\int_{\calU_\b}|\nabla W_{0,M_1,q}|^2\,dx
&\le&C\int_{\Sigma_k}\frac{1}{\sqrt{1-q(\sigma})}\,d\sigma
 \end{eqnarray*}
and the lemma follows at once.
\QED

\section{Existence of $\l^*$}\label{s:localhardy}
 We start with the following local improved Hardy inequality.
\begin{Lemma}\label{lem:loc-hardy}
 Let   $\O$ be a smooth domain and assume that
$\partial\O$ contains a smooth closed submanifold $\Sigma_k$ of
dimension $1\le k\le N-2$. Assume that $p,q$ and $\eta$ satisfy \eqref{eq:weight} and \eqref{eq:min-pq}.
Then there exist constants $\beta_0>0$  and $c>0$
depending only on $\O, \Sig_k,q,\eta$ and $p$ such that for all $\beta\in(0,\beta_0)$
the inequality
$$
\int_{\O_\beta}p|\n
u|^2\,dx-\frac{(N-k)^2}{4}\int_{\O_\beta}q\frac{|u|^2}{\d^{2}}\,dx\geq
c\int_{\O_\beta}\frac{|u|^2}{ \d^{2}|\log(\d)|^{2} }\,dx
$$
holds for all $ u\in H^1_0({\O_\beta})$.
\end{Lemma}
\proof
 We use the notations in Section \ref{s:pn} with $\calU=
\O$ and $\M=\de \O$.\\
Fix $\b_1>0$ small and
 \begin{equation}\label{eq:M2fi}
M_2=-\frac12\,\max_{x\in\ov\O_{\beta_1}}(|h(x)|+ |\n p\cdot \n d |)-1.
 \end{equation}
Since $\frac{p}{q}\in C^1(\ov{\O})$, there exists $C>0$ such that
\be\label{eq:Lippovq}
 \left|\frac{p(x)}{q(x)} - \frac{p(\s(\bar{x}))}{q(\s(\bar{x}))}\right|\leq C\d(x)\quad\forall x\in \O_{\b},
\ee
for small $\b>0$.
Hence by \eqref{eq:min-pq} there exits a constant $C'>0$ such that
\be\label{eq:Lippovq}
p(x)\geq q(x)- C'\d(x)\quad \forall x\in \O_{\b}.
\ee
Consider $ W_{\frac{1}{2},M_2,1}$ (in Lemma \ref{LapFinalExp} with $q\equiv1$).
For  all $\b>0 $ small, we set
\be\label{eq:tiw}
\ti{w}(x)=W_{\frac{1}{2},M_2,1}(x),\quad \forall x\in\O_\b.
\ee
Notice that $\textrm{div}(p\n \ti{w})=p\D \ti{w}+\n p\cdot\n\ti{w}$.
By  Lemma \ref{LapFinalExp}, we have
$$
- \frac{\div (p\n \ti{w})}{\ti{w}}\geq
\frac{(N-k)^2}{4}\,p\d^{-2}+\frac{p}{4}\d^{-2}X_{-2}(\d)
+O({|\log(\d)|\d^{-\frac32}})\,\textrm{ in }\O_\b.
$$
This together with  \eqref{eq:Lippovq} yields
$$
- \frac{\div (p\n \ti{w})}{\ti{w}}\geq
\frac{(N-k)^2}{4}\,q\d^{-2}+\frac{c_0}{4}\d^{-2}X_{-2}(\d)
+O({|\log(\d)|\d^{-\frac32}})\,\textrm{ in }\O_\b,
$$
with $c_0=\min_{\ov{\O_{\b_1}}}p>0$.
Therefore
\be\label{eq:dwow} - \frac{\div (p\n \ti{w})}{\ti{w}}\geq
\frac{(N-k)^2}{4}\,q\d^{-2}+c\,\d^{-2}X_{-2}(\d)\,\textrm{ in
}\O_{\b},
\ee
 for some positive constant $c$ depending only on $\O, \Sig_k,q,\eta$ and $p$.

Let $u\in C^\infty_c(\O_\b)$ and   put
$\psi=\frac{u}{\ti{w}}$. Then one has $|\n
u|^2=|\ti{w}\n\psi|^2+|\psi\n \ti{w}|^2+\n(\psi^2)\cdot \ti{w} \n
\ti{w}$. Therefore $|\n u|^2p=|\ti{w}\n\psi|^2p+p\n
\ti{w}\cdot\n(\ti{w}\psi^2)$. Integrating by parts, we get
$$
\int_{\O_\b}|\n
u|^2p\,dx=\int_{\O_\b}|\ti{w}\n\psi|^2p\,dx+\int_{\O_\b}\left(-
\frac{\textrm{div}(p\n \ti{w})}{\ti{w}}\right)u^2\,dx.
$$
Putting  \eqref{eq:dwow} in the above equality, we get the result.
 \QED

We next prove the following result
\begin{Lemma}\label{lem:Jl1} Let   $\O$ be a smooth bounded domain and assume that
$\partial\O$ contains a smooth closed submanifold $\Sigma_k$ of
dimension $1\le k\le N-2$.  Assume that \eqref{eq:weight} and \eqref{eq:min-pq} hold.
Then there exists $\l^*=\l^*(\O,\Sig_k,p,q,\eta)\in\R$ such
that
$$
\begin{array}{cc}
\displaystyle  \mu_{\l}(\O,\Sigma_k)=\frac{(N-k)^2}{4}, &
\quad\forall
\l\leq\l^*, \vspace{2mm}\\
 \displaystyle \mu_{\l}(\O,\Sigma_k)<\frac{(N-k)^2}{4}, &
\quad\forall \l>\l^*.
\end{array}
$$
\end{Lemma}

\proof   We devide the proof in two steps

\noindent \textbf{Step 1:} We claim that:
\be\label{eq:supmulambda}\sup\limits_{\l\in\R}\mu_\l(\O,\Sigma_k)\leq
\frac{(N-k)^2}{4}. \ee Indeed, we know that
$\nu_0(\R^N_+,\R^k)=\frac{(N-k)^2}{4}$, see  \cite{FTT} for instance. Given
$\tau>0$, we let $u_\tau\in C^\infty_c(\R^N_+)$ be such that
\begin{equation}\label{eq:estutau}
\int_{\R^N_+}|\n
u_\tau|^2\,d y\leq\left(\frac{(N-k)^2}{4}+\tau\right)\int_{\R^N_+}|\tilde
y|^{-2}u_\tau^2\,d y.
\end{equation}
By \eqref{eq:min-pq}, we can let $\sigma_0\in\Sigma_k$ be such that
$$
q(\sigma_0)=p(\s_0).
$$
Now, given $r>0$, we let  $\rho_r>0$ such that for all $ x\in
B(\sigma_0,\rho_r)\cap \Omega $
\begin{equation}\label{eq:estq}
p(x)\le (1+r)q(\sigma_0),\quad   q(x)\ge (1-r)q(\sigma_0)\quad\textrm{ and }\quad \eta(x)\le r.
\end{equation}
We choose Fermi coordinates near $\sigma_0\in\Sigma_k$ given by the map $ F^{\sigma_0}_{\partial\O}$ (as in Section \ref{s:pn}) and we choose
$\e_0>0$ small such that, for all $\e\in(0,\e_0) $,
$$
\Lambda_{\e,\rho,r,\tau}:=F^{\sigma_0}_{\partial\O}(\e\,{\rm
Supp(u_\tau)})\subset\,B(\sigma_0,\rho_r)\cap \Omega
$$
and we define the following test function
$$
v(x)=\e^{\frac{2-N}{2}}u_\tau\left(\e^{-1}(F^{\sigma_0}_{\partial\O})^{-1}(x)\right),
\quad x\in \Lambda_{\e,\rho,r,\tau}.
$$
Clearly, for every $\e\in(0,\e_0)$, we have that $v\in
C^\infty_c(\O)$  and thus by a change of variable, \eqref{eq:estq}
and Lemma \ref{lemddelta}, we have~
\begin{eqnarray*}
\mu_\l(\O,\Sigma_k)&\leq&\frac{\displaystyle \int_{\O}p|\n v|^2\,dx
+\l\int_{\O}\d^{-2}\eta v^2\,dx}{\displaystyle
\int_{\O}q(x)\,\d^{-2}\,v^2\,dx}\\
&\leq&\frac{\displaystyle (1+r)\int_{\Lambda_{\e,\rho,r,\tau}}|\n
v|^2\,dx
}{(1-r)\,\displaystyle
\int_{\Lambda_{\e,\rho,r,\tau}}\d^{-2}\,v^2\,dx}+\frac{r|\l|}{(1-r)q(\sigma_0) } \\
&\leq&\frac{\displaystyle (1+r)\int_{\Lambda_{\e,\rho,r,\tau}}|\n
v|^2\,dx
}{(1-c r)\,\displaystyle
\int_{\Lambda_{\e,\rho,r,\tau}}\ti{\d}^{-2}\,v^2\,dx}+\frac{r|\l|}{(1-r)q(\sigma_0) } \\
&\leq&\frac{(1+r)\e^{2-N}\displaystyle
\int_{\R^N_+}\e^{-2}(g^\e)^{ij}\partial_i
u_\tau\partial_ju_\tau|\,\sqrt{|g^\e|}(y)\,dy
}{(1-cr)\,\displaystyle
\int_{\R^N_+}\e^{2-N}\,|\e\tilde y|^{-2}\,u_\tau^2\,\sqrt{|g^\e|}(\tilde y)\,d y}+\frac{cr}{1-r } ,\\
\end{eqnarray*}
where $g^\e$ is the scaled metric with components $g^\e_{\a\b}(y)=\e^{-2}\la\de_\a F^{\s_0}_{\de\O}(\e y), \de_\b F^{\s_0}_{\de\O}(\e y)\ra$
for $\a,\b=1,\dots,N$
and where we have used the fact that $\tilde{\d}(F^{\s_0}_{\de\O}(\e y))=|\e\tilde  y|^2$ for 
every $\tilde y$ in the support of $u_\t$.
Since the scaled metric $g^\e$ expands a $g^\e=I+O(\e)$ on the support of $u_\t$, we deduce that
\begin{eqnarray*}
 \mu_\l(\O,\Sigma_k) &\le&  \frac{1+r}{1-c r}\,\frac{1+c\e}{1-c\e}\,\, \frac{\displaystyle
\int_{\R^N_+}|\nabla u_\tau|^2\,d y }{\displaystyle
\int_{\R^N_+}|\tilde y|^{-2}\,u_\tau^2\,d y}+\frac{cr}{1-r}   ,
\end{eqnarray*}
where $c$  is a positive constant depending only on $\O,p,q,\eta$ and $\Sig_k$. Hence by \eqref{eq:estutau} we conclude
\begin{eqnarray*}
 \mu_\l(\O,\Sigma_k)
&\le& \frac{1+r}{1-c r}\,\frac{1+c\e}{1-c\e}\, \left( \frac{(N-k)^2}{4}+\tau
\right)+ \frac{cr}{1-r}  .
\end{eqnarray*}
Taking the limit in $\e$, then in $r$ and then in $\tau$, the claim follows.\\
 \textbf{Step 2:} We claim that there exists $\ti{\l}\in\R$ such that
 $\mu_{\ti{\l}}(\O,\Sig_k)\geq\frac{(N-k)^2}{4}$.\\
Thanks to Lemma \ref{lem:loc-hardy}, the proof uses a standard argument of cut-off function
 and integration by parts (see \cite{BM}) and we can obtain
$$
\int_{\O}\d^{-2}u^2 q\,dx\leq \int_{\O}|\n u|^2 p\,dx+C\int_{\O}\d^{-2}u^2 \eta \,dx\quad\forall u\in C^\infty_c(\O),
$$
for some constant $C>0$.  We skip the details. The claim now follows by choosing $\ti{\l}=-C$\\

Finally, noticing that  $\mu_\l(\O,\Sig_k)$ is
  decreasing in $\l$,   we can set
\be\label{eq:lsdef} \l^*:=\sup\left\{{\l\in\R}\,:\, \mu_\l(\O,\Sig_k)=
{\frac{(N-k)^2}{4}}\right\}
 \ee
 so that $\mu_\l(\O,\Sig_k)<\frac{(N-k)^2}{4}$ for all $\l>\l^*$.
\QED

\section{Non-existence result}\label{s:ne}

\begin{Lemma}\label{lem:Opm}
Let $\O$ be a smooth bounded domain of $\R^N$, $N\geq 3$, and let $\Sigma_k$ be a
smooth closed submanifold of  $\partial\O$ of dimension $k$ with $1\le k\le
N-2$. Then, there exist bounded  smooth domains $\O^\pm$ such that $\O^+\subset \O\subset\O^-$
and
$$
\de{\O^+}\cap \de\O=\de{\O^-}\cap \de\O = \Sigma_k.
$$
\end{Lemma}
\proof
Consider the  maps
$$
x\mapsto g^\pm(x):=d_{\partial\O}(x)\pm\frac12\,\d^2(x),
$$
where $d_{\de\O}$ is the distance function to $\de\O$.
For some $\b_1>0$ small, $g^\pm$ are smooth in $\O_{\b_1}$ and since $|\n g^\pm|\geq C>0$ on $\Sig_k$, by the implicit function theorem, the sets
$$
\{x\in \O_{\b}\,:\,g^\pm=0 \}
$$
are smooth $(N-1)$-dimensional submanifolds of $\R^N$, for some $\b>0$ small. In addition, by construction, they can be taken to be part of the boundaries of
smooth bounded domains  $\O^\pm $ with  $\O^+\subset \O\subset\O^-$ and such that
$$
\de{\O^+}\cap \de \O=\de{\O^-}\cap \de\O = \Sigma_k.
$$
The prove then follows at once.
 \QED

Now, we prove the following non-existence result.
\begin{Theorem}\label{th:ne}
Let $\O$ be a smooth bounded domain of $\R^N$ and let $\Sigma_k$ be a
smooth closed submanifold of  $\partial\O$ of dimension $k$ with  $1\le k\le
N-2$ and let $\l\geq0$.  Assume that $p,q$ and $\eta$ satisfy \eqref{eq:weight} and \eqref{eq:min-pq}. Suppose that  $u\in H^1_0(\O)\cap C(\O)$ is
a non-negative function
 satisfying
\be\label{eq:ustf}
-\div(p \n u)-\frac{(N-k)^2}{4}q\d^{-2}u\geq-\l \eta \d^{-2} u \quad\textrm{in }\O.
\ee
If $\int_{\Sigma_k}\frac1{\sqrt{1-p(\s)/q(\s)}}d\s=+\infty$ then
$u\equiv0$.
\end{Theorem}
\proof
We first assume that $p\equiv1$.
 Let $\O^+$ be the set given by  Lemma \ref{lem:Opm}. We
will use the notations in Section \ref{s:pn} with $\calU=
\O^+$ and $\M=\de \O^+$. For  $\b>0$ small we define
$$\O^+_{\b} := \{ x \in \O^+: \quad
{\d}(x)<\b \}.$$
 We suppose by contradiction  that $u $ does not vanish identically  near $\Sigma_k$  and  satisfies
 \eqref{eq:ustf} so that $u>0$ in  $\O_{\b}$ by  the maximum principle,  for some $\b>0$ small.\\
Consider the subsolution  $V_\e$ defined in  Lemma \ref{le:lowerbound} which satisfies
\be\label{eq:lwaneg}
\calL_\l\,V_\e\leq 0\quad\textrm{ in
 }\O^+_{\b},\quad\forall \e\in(0,1).
\ee
Notice that $\ov{\de\O^+_{\b}\cap
\O^+}\subset \O$ thus, for  $\b>0$ small,  we can choose $R>0$ (independent on $\e$) so
that
$$
R\,V_\e\leq R\,V_0\leq u\quad\textrm{ on }
\ov{\de\O^+_{\b}\cap \O^+ }
\quad\forall \e\in(0,1).
$$
Again by Lemma \ref{le:lowerbound}, setting $v_\e=R\, {V_\e}-u$, it
turns out that $v^+_\e=\max(v_\e,0)\in
H^1_0(\O^+_{\b})$ because $V_\e=0$ on $\de
\O^+_{\b}\setminus
\ov{\de\O^+_{\b}\cap
\O^+}$. Moreover by \eqref{eq:ustf} and
\eqref{eq:lwaneg},
$$
\calL_\l\,v_\e\leq 0\quad\textrm{ in
}\O^+_{\b},\quad\forall \e\in(0,1).
$$
 Multiplying the above inequality by $v^+_\e$ and integrating by parts
 yields
$$
\int_{\O^+_{\b}}|\n
v^+_\e|^2\,dx-\frac{(N-k)^2}{4}\int_{\O^+_{\b}}\d^{-2}q|v^+_\e|^2\,dx+
\l\int_{\O^+_{\b}}\eta \d^{-2}|v^+_\e|^2\,dx
\leq0.
$$
But then Lemma \ref{lem:loc-hardy} implies that $v^+_\e=0$ in
$\O^+_{\b}$ provided $\b$ small enough because $|\eta|\leq C\d$ near $\Sig_k$. Therefore $u\geq R\, {V_\e}$ for every $\e\in(0,1)$. In
particular $u\geq R\,V_0$. Hence we obtain from Lemma \ref{le:lowerbound} that
$$
\infty>\int_{\O^+_{\b}}\frac{u^2}{\d^{2}}\geq R^2 \int_{\O^+_{\b}}\frac{V_0^2}{\d^{2}}\geq\int_{\Sigma_k}\frac1{\sqrt{1-q(\s)}}d\s
$$
which leads to a contradiction. We deduce that $u\equiv0$ in $\O^+_{\b} $. Thus by
the maximum principle $u\equiv0$ in  $\O$.\\
For the general case $p\neq 1$, we argue as in \cite{BMS} by setting
\be\label{eq:transf}
\ti{u}=\sqrt{p} u.
\ee
 This function satisfies
$$
-\D \ti{u}-\frac{(N-k)^2}{4}\frac{q}{p}\d^{-2}\ti{u}\geq-\l \frac{\eta}{p} \d^{-2}\ti{ u} +\left(-\frac{\D p}{2 p } +\frac{|\n p|^2}{4 p^2 } \right) \ti{u}\quad\textrm{in }\O.
$$
Hence since $p\in C^2(\ov{\O})$ and $p>0$ in $ \ov{\O}$, we get the same conclusions as in the case $p\equiv 1$ and $q$ replaced by $q/p$.
 \QED
%
\section{Existence of minimizers for $\m_{\l}(\Omega,\Sigma_k) $}

\begin{Theorem}\label{th:exitslesls}
Let $\O$ be a smooth bounded domain of $\R^N$ and let $\Sigma_k$ be a
smooth closed submanifold of  $\partial\O$ of dimension $k$ with $1\le k\le
N-2$.  Assume that $p,q$ and $\eta$ satisfy \eqref{eq:weight} and \eqref{eq:min-pq}. Then  $\m_{\l}(\Omega,\Sigma_k)$ is achieved for every $\l<\l^*$.
\end{Theorem}
\proof
The proof follows the same argument of \cite{BM} by taking into account the fact that $\eta=0$ on $\Sig_k$ so we skip it.
\QED

Next, we prove the existence of minimizers in the critical case $\l=\l_*$.
\begin{Theorem}\label{th:exits-crit}
Let $\O$ be a smooth bounded domain of $\R^N$ and let $\Sigma_k$ be a
smooth closed submanifold of  $\partial\O$ of dimension $k$ with  $1\le k\le
N-2$.  Assume that $p,q$ and $\eta$ satisfy \eqref{eq:weight} and \eqref{eq:min-pq}. If $\displaystyle \int_{\Sigma_k}\frac1{\sqrt{1-p(\s)/q(\s)}}d\s<\infty$ then
$\m_{\l^*}=\m_{\l^*}(\O,\Sig_k)$ is achieved.
\end{Theorem}
\proof
We first consider the case $p\equiv 1$.\\
Let $\l_n$ be a sequence of real numbers decreasing to $\l^*$. By Theorem \ref{th:exitslesls}, there exits $u_n$
minimizers for $\mu_{\l_n}=\m_{\l_n}(\Omega,\Sigma_k)$ so that
\be\label{eq:u_n}
-\D u_n-\mu_{\l_n}\d^{-2}q u_n= -\l_n  \d^{-2 }\eta u_n \quad\textrm{ in }\O.
\ee
We may assume that $u_n\geq 0$ in $\O$. We may also assume  that $\|\n u_n\|_{L^2(\O)}=1$. Hence $u_n \rightharpoonup u$ in $H^1_0(\O)$
and $u_n\to u$ in $L^2(\O)$ and pointwise.
Let $\O^-\supset\O$ be the set given by  Lemma \ref{lem:Opm}. We
will use the notations in Section \ref{s:pn} with $\calU=
\O^-$ and $\M=\de \O^-$. It will be understood that $q$ is extended to a function in  $C^2(\ov{\O^- })$.
 For  $\b>0$ small we define
$$\O^-_{\b} := \{ x \in \O^-: \quad
{\d}(x)<\b \}.$$
We have that
$$
\D u_n+b_n(x)\, u_n=0\quad\textrm{ in }\O,
$$
with $|b_n|\leq C$ in $\ov{\O\setminus \ov{\O^-_{\frac\b2}}}$ for all integer $n$. Thus by standard elliptic regularity theory,
\be\label{eq:unleC}
u_n\leq C \quad \quad\textrm{ in }\ov{\O\setminus \ov{\O^-_{\frac{\b}{2}}}}.
\ee
We consider the supersolution $U$ in  Lemma \ref{le:upperbound}. We shall show that there exits a constant $C>0$ such that for all $n\in\N$
\be\label{eq:unleCV12}
u_n\leq C U \quad \textrm{ in }\ov{\O^-_\b}.
\ee
Notice that $\ov{\O\cap\de\O^-_\b}\subset \O^-$  thus by \eqref{eq:unleC}, we can choose $C>0$ so
that for any $n$
$$
 u_n\leq C\, U\quad\textrm{ on }
\ov{\O\cap\de\O^-_\b}.
$$
Again by Lemma \ref{le:upperbound}, setting $v_n=u_n-C\, U$, it
turns out that $v^+_n=\max(v_n,0)\in
H^1_0(\O^-_{\b})$ because $u_n=0$ on $\de\O\cap \O^-_\b$.
Hence we have
$$
\calL_{\l_n}\,v_n\leq -C(\mu_{\l^*}-\mu_n)q{U}-C(\l^*-\l_n)\eta {U}\leq 0 \quad\textrm{ in
}\O^-_{\b}\cap\O .
$$
 Multiplying the above inequality by $v^+_n$ and integrating by parts
 yields
$$
\int_{\O^-_{\b}}|\n
v^+_n|^2\,dx-\mu_{\l_n}\int_{\O^-_{\b}}\d^{-2}q|v^+_n|^2\,dx+
\l_n\int_{\O^-_{\b}}\eta\d^{-2} |v^+_n|^2\,dx
\leq0.
$$
Hence Lemma \ref{lem:loc-hardy} implies that
$$
C \int_{\O^-_{\b}}\d^{-2}X_{-2} |v^+_n|^2\,dx+\l_n\int_{\O^-_{\b}}\eta \d^{-2} |v^+_n|^2\,dx\leq0.
$$
Since $\l_n$ is bounded, we can choose $\b>0$ small (independent of
$n$) such that $v^+_n\equiv0$ on $\O^-_\b$ (recall that $|\eta|\leq
C\d$).
 Thus we obtain \eqref{eq:unleCV12}. \\
Now since $u_n\to u$ in $L^2(\O)$, we get by the dominated convergence theorem and \eqref{eq:unleCV12}, that
$$
\d^{-1} u_n\to  \d^{-1} u\quad \textrm{ in }L^2(\O).
$$
Since $u_n$ satisfies
$$
1=\int_{\O}|\n u_n|^2=\mu_{\l_n}\int_{\O}\d^{-2} q u_n^2+ {\l_n}\int_{\O}\d^{-2} \eta u_n^2,
$$
taking the limit, we have $1= \mu_{\l^*}\int_{\O}\d^{-2} q u^2+ {\l^*}\int_{\O}\d^{-2} \eta u^2$. Hence $u\neq0$ and it is a minimizer for  $\mu_{\l^*}=\frac{(N-k)^2}{4}$.\\
For the general case $p\neq 1$, we can use the same transformation as in \eqref{eq:transf}. So \eqref{eq:unleCV12} holds and the same argument as a above carries over.
\QED
\section{Proof of Theorem \ref{th:mulpqe} and Theorem \ref{th:crit}}
\textit{Proof of Theorem \ref{th:mulpqe}:} Combining Lemma \ref{lem:Jl1} and Theorem \ref{th:exitslesls},
it remains only to check the case $\l<\l^*$. But this is an easy consequence of the definition of $\l^*$ and of $\mu_{\l}(\O,\Sig_k)$, see [\cite{BM}, Section 3].\QED
\bigskip
\noindent
\textit{Proof of Theorem \ref{th:crit}:}
Existence is proved in Theorem \ref{th:exits-crit} for $I_k<\infty$.  Since the absolute value of
any minimizer for $\mu_{\l}(\O,\Sig_k) $ is also a minimizer, we can apply Theorem \ref{th:ne}
to infer that $\mu_{\l^*}(\O,\Sig_k) $ is never achieved as soon as $I_k=\infty$.\QED

%


\begin{center}\textbf{ Acknowledgments} \end{center}
This work started when the first author was visiting CMM,
Universidad de Chile. He is grateful for their kind hospitality. M.
M. Fall is supported by the Alexander-von-Humboldt Foundation. F.
Mahmoudi is supported by the Fondecyt proyect n: 1100164 and Fondo
Basal CMM.

\end{document}